\documentclass[10pt, a4paper]{article}

\usepackage{amsmath,  amsfonts}

\usepackage[T1]{fontenc}

\newtheorem{theorem}{\quad Theorem}[section]

\newcommand{\be} {\begin{equation}}
\newcommand{\ee} {\end{equation}}

\title {Uniform bound for the volume of the solutions to Liouville type equations on the annulus.}

\date{}

 \author{Samy Skander Bahoura\footnote {e-mails: samybahoura@yahoo.fr, samybahoura@gmail.com} \\ 
 {\small Equipe d'Analyse Complexe et G\'eom\'etrie.}\\  
  {\small Universit\'e Pierre et Marie Curie, 75005 Paris, France.}} 

\begin{document}

\maketitle
\begin{abstract}

We consider variational problems with regular H\"olderian weight or with weight and boundary singularity and, Dirichlet condition. We prove the boundedness of the volume of the solutions to these equations on the annulus.

\end{abstract}



\section{Introduction and Main Results} 

We set $ \Delta = \partial_{11} + \partial_{22} $  on the annulus $ \Omega = C(1,1/2,0) $ of $ {\mathbb R}^2 $ of radii $ 1 $ and $ 1/2 $ centered at the origin.

\smallskip

We consider the following equation:

$$ (P)   \left \{ \begin {split} 
      -\Delta u & = (1+|x-x_0|^{2\beta}) V e^{u} \,\, &\text{in} \,\, & \Omega  \subset {\mathbb R}^2, \\
                  u & = 0  \,\,             & \text{in} \,\,    &\partial \Omega.              
\end {split}\right.
$$

Here: $ C(1) $ the unit circle and $ C(1/2) $ the circle of radius $ 1/2 $ centered at the origin.

$$ \beta \in (0, 1/2), \,\, x_0 \in C(1/2), $$
and,

$$ u \in W_0^{1,1}(\Omega), \,\, e^u \in L^1({\Omega}) \,\, {\rm and} \,\,  0 < a \leq V \leq b. $$

This is an equation with regular H\"olderian weight not Lipschitz in $ x_0 $ but have a weak derivative.

We have in [10] a nice formulation of this problem $ (P) $ in the sense of the distributions. This Problem arises in various geometrical and physical situations, see for example [1, 5, 23, 27]. The above equation was studied by many authors, with and without boundary conditions, also for Riemannian surfaces, see [1-27],  where one can find some existence and compactness results. In [9] we have the following important Theorem,

{\bf Theorem A}{\it (Brezis-Merle [9])}.{\it For $ u $ and $ V $ two functions relative to $ (P) $ with,
$$ 0 < a \leq V \leq b < + \infty $$
then it holds, for all compact set $ K $ of $ \Omega $:

$$ \sup_K u \leq c, $$

with $ c $ depending on $ a, b, \beta, x_0, K $ and $ \Omega $.}

We deduce from Theorem A and from the elliptic estimates that, $ u $ is uniformly bounded in $ C_{loc}^2(\Omega) $.

In this paper we try to prove that we have on all $ \Omega $ the boundedness of the volume of the solutions of $ (P) $ if we add the assumption that $ V $ is uniformly Lipschitz with particular Lispchitz number.

Here we have:

\begin{theorem} Assume that $ u $ is a solution of $ (P) $ relative to $ V $ with the following conditions:

$$  x_0 \in C(1/2) \subset \partial \Omega, \,\, \beta \in (0,1/2), $$
and,
$$ 0 < a \leq V \leq b,\,\, ||\nabla V||_{L^{\infty}} \leq A=\dfrac{a}{2(1+2^{2\beta})}, $$
we have,

$$  \int_{\Omega} e^u \leq c(a, b, \beta, x_0, \Omega). $$

\end{theorem} 

\smallskip

A consequence of this theorem is a compactness result of the solutions to this Liouville type equation, see [4].

\smallskip

We have the same result if we consider Liouville equation with boundary singularity on the annulus $ \Omega=C(1,1/2,0) $ for particular Lipschitz number of the potential $ V $:

We consider the following equation in the sense of distributions:

$$ (P_{\beta})   \left \{ \begin {split} 
      -\Delta u & = |x-x_0|^{2\beta} V e^{u} \,\, &\text{in} \,\, & \Omega  \subset {\mathbb R}^2, \\
                  u & = 0  \,\,             & \text{in} \,\,    &\partial \Omega.              
\end {split}\right.
$$

Here: $ C(1) $ the unit circle and $ C(1/2) $ the circle of radius $ 1/2 $ centered at the origin.

$$ \beta \in (-1/2, +\infty), \,\, x_0 \in C(1/2), $$
and,

$$ u \in W_0^{1,1}(\Omega), \,\, |x-x_0|^{2\beta} e^u \in L^1({\Omega}) \,\, {\rm and} \,\,  0 < a \leq V \leq b. $$

This is an equation with boundary singularity and Dirichlet condition on the annulus.

Here we have:

\begin{theorem} Assume that $ u $ is a solution of $ (P_{\beta}) $ relative to $ V $ with the following conditions:

$$  x_0 \in C(1/2) \subset \partial \Omega, \,\, \beta \in (-1/2, +\infty), $$

and,

$$ 0 < a \leq V \leq b,\,\, ||\nabla V||_{L^{\infty}} \leq A=\dfrac{(\beta+1) a}{2}, $$

we have,

$$  \int_{\Omega} |x-x_0|^{2\beta} e^u \leq c(a, b, \beta, x_0, \Omega). $$

\end{theorem} 

For the previous theorem, the condition $ \beta >-1/2 $ is in relation with the regularity of the solutions to apply the Pohozaev identity.

A consequence of the previous theorem is a compactness result for the solutions to a Liouville type equation with boundary singularity on the annulus, see [3].

\section{Proof of the Theorems:}

\smallskip

{\it Proof of the theorems:}

\smallskip

By corollary 1 of the paper of Brezis-Merle, we have: $ e^{ku} \in L^1(\Omega) $ for all $ k >2 $ and the elliptic estimates and the Sobolev embedding imply that: $ u \in W^{2,k}(\Omega) \cap C^{1,\epsilon}(\bar \Omega), \epsilon >0 $. By the maximum principle $ u\geq 0 $.

\smallskip

Step 1: By using the first eigenvalue and the first eigenfunction, with Dirichlet boundary condition, the volume is locally uniformly bounded, and thus the solutions are locally uniformly bounded by Brezis-Merle result. The solutions $ u >0 $ are locally uniformly bounded in $ C^{1,\epsilon}(\Omega) $ for $ \epsilon $ small.

\smallskip

Step 2: Let's consider $ C_1=C(1,3/4,0) $ and $ C_2=C(3/4,1/2,0) $ the two annulus wich are the neighborhood of the two components of the boundary.

We multiply the equation by $ (x-x_0)\cdot \nabla u $ on $ C_1 $ and $ C_2 $ and use the Pohozaev identity. We use the fact that $ u $ is uniformly bounded around the circle $ C(3/4) $. We obtain:

1) We have on $ C_1 $:

$$ \int_{C_1} (\Delta u)[(x-x_0)\cdot \nabla u] dx=\int_{C_1} -[(1+|x-x_0|^{2\beta}) V (x-x_0)\cdot \nabla (e^u)]dx, $$

Thus,

$$
 \int_{\partial C_1} [(x-x_0)\cdot \nabla u ] (\nabla u \cdot \nu)-\dfrac{1}{2}[(x-x_0)\cdot \nu] |\nabla u|^2= $$
 
$$ = \int_{C_1} (2+2(\beta + 1)|x-x_0|^{2\beta}) Ve^u dx + \int_{C_1} (1+|x-x_0|^{2\beta}) (x-x_0)\cdot \nabla V e^u dx+  $$

$$ - \int_{\partial C_1} (1+|x-x_0|^{2\beta}) [(x-x_0)\cdot \nu ] V e^u d\sigma $$

We can write, ($ u= 0 $ on $ C(1) $):

$$ \int_{C(1)} \frac{1}{2} [(x-x_0)\cdot \nu ] (\partial_{\nu} u)^2 d\sigma + O(1)= $$

$$ =\int_{C_1} (2+2(\beta + 1)|x-x_0|^{2\beta}) Ve^u dx + \int_{C_1} (1+|x-x_0|^{2\beta}) (x-x_0)\cdot \nabla V e^u dx+ O(1) = $$

$$ = C\int_{C(1)} \partial_{\nu} u d\sigma + O(1), $$

with $ C>0 $ not depends on $ u $.

Be cause $ \nu = x, ||x||=1, ||x_0||=1/2 $ and then by Cauchy-Schwarz, $ (x-x_0)\cdot x=||x||^2-x_0 \cdot x \geq 1/2 $, we obtain:

$$ \int_{C(1)} (\partial_{\nu} u )^2 d\sigma \leq 4 C\int_{C(1)} \partial_{\nu} u d\sigma +O(1), $$

and, by the Cauchy-Schwarz inequality applied to the integral of the right hand side we obtain:

$$ \int_{C(1)} (\partial_{\nu} u)^2 d\sigma = O(1),\,\, {\rm and } \,\, \int_{C(1)} (\partial_{\nu} u) d\sigma =O(1). $$

2) We have on $ C_2 $, we use again, the uniform boundedness of $ u $ in $ C^1 $ norm around $ C(3/4) $:

$$ \int_{C_2} (\Delta u)[(x-x_0)\cdot \nabla u] dx=\int_{C_2} -[(1+|x-x_0|^{2\beta}) V (x-x_0)\cdot \nabla (e^u)]dx, $$

Thus,

$$
 \int_{\partial C_2} [(x-x_0)\cdot \nabla u ] (\nabla u \cdot \nu)-\dfrac{1}{2}[(x-x_0)\cdot \nu] |\nabla u|^2= $$
 
$$ = \int_{C_2} (2+2(\beta + 1)|x-x_0|^{2\beta}) Ve^u dx + \int_{C_2} (1+|x-x_0|^{2\beta}) (x-x_0)\cdot \nabla V e^u dx+  $$

$$ - \int_{\partial C_2} (1+|x-x_0|^{2\beta}) [(x-x_0)\cdot \nu ] V e^u d\sigma $$

We can write,($ u=0 $ on $ C(1/2) $):

$$ \int_{C(1/2)} \frac{1}{2} [(x-x_0)\cdot \nu ] (\partial_{\nu} u)^2 d\sigma + O(1)= $$

$$ = \int_{C_2} (2+2(\beta +1)|x-x_0|^{2\beta}) Ve^u dx + \int_{C_2} (1+|x-x_0|^{2\beta}) (x-x_0)\cdot \nabla V e^u dx+ O(1), $$

But here, $ \nu=-2x$ and $ (x-x_0)\cdot \nu= -(x-x_0)\cdot x \leq 0 $ and thus:

$$ \int_{C_2} (2+2(\beta +1)|x-x_0|^{2\beta}) Ve^u dx + \int_{C_2} (1+|x-x_0|^{2\beta}) (x-x_0)\cdot \nabla V e^u dx + $$

$$ + \int_{C(1/2)} \frac{1}{2}[-(x-x_0)\cdot \nu ] (\partial_{\nu} u)^2 d\sigma= O(1), $$

If we choose:

\be \dfrac{|(x-x_0)\cdot \nabla V|}{V} \leq \dfrac{1}{2} \inf_{x \in \bar \Omega } \dfrac{(2+2(\beta+1)|x-x_0|^{2\beta})}{1+|x-x_0|^{2\beta}}, \ee

We obtain:

$$ \int_{C_2} (2+2(\beta +1)|x-x_0|^{2\beta}) Ve^u dx + \int_{C_2} (1+|x-x_0|^{2\beta}) (x-x_0)\cdot \nabla V e^u dx \geq $$

$$ \geq \frac{1}{2} \int_{C_2} (2+2(\beta +1)|x-x_0|^{2\beta}) Ve^u dx \geq 0 $$

thus,

$$ \int_{C_2} (2+2(\beta +1)|x-x_0|^{2\beta}) Ve^u dx =O(1), $$

(This condition is true for the second theorem for $ x \in \Omega $, we have: $ \dfrac{|(x-x_0)\cdot \nabla V|}{V} \leq \dfrac{1}{2} \inf_{x \in \Omega } \dfrac{2(\beta+1)|x-x_0|^{2\beta}}{|x-x_0|^{2\beta}} =\beta+1 $, this is true if $ A \leq \dfrac{a(\beta + 1)}{2} $.)

and thus,

$$ \int_{C(1/2)} \partial_{\nu} u d\sigma = O(1), $$

Thus, if: $ 2 \dfrac{diam (\Omega) ||\nabla V||_{\infty}}{a} \leq \dfrac{2}{1+[diam (\Omega)]^{2\beta}}  \Rightarrow $ if $2 A\leq \dfrac{2a}{2(1+[diam(\Omega)]^{2\beta})} $, we obtain a uniform bound for $\int_{C(1/2)} \partial_{\nu} u d\sigma $. But $ diam(\Omega) = diameter (\Omega) =2 $, because it is the annulus.

If we use  1) and 2), we obtain: if $  A\leq \dfrac{a}{2(1+2^{2\beta})} $ for the first theorem, then:

$$\int_{C(1)} (\partial_{\nu} u) d\sigma =O(1),\,\, {\rm and} \,\, \int_{C(1/2)} \partial_{\nu} u d\sigma = O(1),$$

and, thus:

$$ \int_{\Omega } (1+|x-x_0|^{2\beta}) Ve^u dx=\int_{\partial \Omega} \partial_{\nu} u d\sigma =O(1). $$

We have the same conclusion for the second theorem if $ A=\dfrac{a(\beta +1)}{2}$.

\end{document}